\documentclass[12pt]{article}
\newcommand\norm[1]{\left\lVert#1\right\rVert}

\usepackage{amsthm}
\usepackage{amsmath}
\usepackage{enumerate}
\usepackage{amssymb}
\usepackage{latexsym}
\usepackage{amsfonts}
\usepackage{color}
\usepackage{mathrsfs}
\usepackage{epsfig}
\newtheorem{theorem}{\bf Theorem}[section]
\newtheorem{proposition}[theorem]{\bf Proposition}
\newtheorem{definition}[theorem]{\bf Definition}

\newtheorem{remark}[theorem]{\bf Remark}
\newtheorem{lemma}[theorem]{\bf Lemma}

\newsavebox{\savepar}

\begin{document}
	\title{Existence results of two mixed boundary value elliptic PDEs in $\mathbb{R}^n$}
	\author{ Akasmika Panda\footnote{Corresponding
			author: akasmika44@gmail.com}~ \& Debajyoti Choudhuri \\
\small{Department of Mathematics, National Institute of Technology Rourkela}\\
	\small{Emails: akasmika44@gmail.com, dc.iit12@gmail.com}
	}
	\date{}
	\maketitle
	\begin{abstract}
		\noindent We study the existence of a solution to the mixed boundary value problem for Helmholtz and Poisson type equations in a bounded Lipschitz domain $\Omega\subset\mathbb{R}^N$ and in $\mathbb{R}^N\setminus\Omega$ for $N\geq3$. The boundary $\partial\Omega$ of $\Omega$ is the decomposition of $\Gamma_1,\Gamma_2\subset\partial\Omega$ such that  $\partial\Omega=\Gamma=\overline{\Gamma}_1\cup\Gamma_2=\Gamma_1\cup\overline{\Gamma}_2$ and $\Gamma_1\cap\Gamma_2=\emptyset$. We have shown that if the Neumann data $f_2\in H^{-\frac{1}{2}}(\Gamma_2)$ and the Dirichlet data  $f_1\in H^{\frac{1}{2}}(\Gamma_1)$ then the Helmholtz problem with mixed boundary data admits a unique solution. We have also shown the existence of a weak solution to a mixed boundary value problem governed by the Poisson equation with a measure data and the Dirichlet, Neumann data belongs to  $f_1\in H^{\frac{1}{2}}(\Gamma_1)$,  $f_2\in H^{-\frac{1}{2}}(\Gamma_2)$ respectively.
			\\
		{\bf keywords}:~Mixed boundary value problem, Sobolev space, Newton potential, Boundary integral operator, Layer potentials, Radon measure.
	\end{abstract}
	\section{Introduction}
	The Poisson problem with mixed Dirichlet-Neumann boundary conditions deals with conductivity, heat transfer, metallurgical melting, wave phenomena, elasticity and electrostatics in mathematical physics and engineering. The detailed applications can be found in \cite{Dauge}, \cite{Fabrikant}, \cite{Jochmann}, \cite{Lagnese}, \cite{Mazya}, \cite{Rempel}, \cite{Simanca}, \cite{Sneddon}, \cite{Wendland} and the references therein. A common problem of interest found in the literature is the following mixed boundary value (MBVP).
	\begin{eqnarray}\label{m1}
	\begin{split}
	L u&= h~\text{in}~\Omega,\\
	u&= f~\text{on}~\Gamma_1,\\
	Mu&= g~\text{on}~\Gamma_2,
	\end{split}
	\end{eqnarray}
	where, $\Omega$ is a bounded Lipschitz domain in $\mathbb{R}^N$ for $N\geq3$. The boundary of $\Omega$, which will be denoted by $\Gamma$, is the disjoint union of $\Gamma_1$ and $\Gamma_2$ which are subsets of $\Gamma$ such that  $\overline{\Gamma}_1\cup\Gamma_2=\Gamma_1\cup\overline{\Gamma}_2=\Gamma$ and $\Gamma_1\cap\Gamma_2=\emptyset$. Further, $L$ is a second order elliptic operator, $M$ is a general first order oblique differential operator on $\Gamma_2$.\\
	Lieberman \cite{Lieberman 1,Lieberman 2} considered the problem $\eqref{m1}$ and proved the existence and H\"{o}lder continuity of classical solutions with smooth data. The techniques used in the corresponding Dirichlet problem ($\Gamma_2=\emptyset$) and oblique derivative problem ($\Gamma=\Gamma_2$) of $\eqref{m1}$ are helpful to show the existence of solutions to the mixed boundary value problem. It is worth to mention the work due to Azzam and Kreyszig \cite{Azzam}, as they have provided the regularity result for MBVP in a plane domain with corners, where the Dirichlet data belongs to $C^{2,\alpha}(\Gamma\setminus\{0\})$ and the remaining boundary data is in $C^{1,\alpha}(\Gamma\setminus\{0\})$. The work due to Sykes \cite{Sykes} deals with the boundary regularity of problem $\eqref{m1}$ with Dirichlet and Neumann boundary conditions where $L$ is the Laplacian operator, $h=0$ in $\Omega$,
	 $g\in L^p(\Gamma_2)$, $f\in W^{1,p}(\Gamma_1)$ for $1<p\leq2$ and the angle  between $\Gamma_1$, $\Gamma_2$ should be strictly less than $\pi$ in the interface. Sykes \cite{Sykes} drew motivation from Brown \cite{Brown}. who considered the two boundary data as $f\in H^1(\Gamma_1)$ and $g\in L^2(\Gamma_2)$. \\
	Not much of literature is found for MBVP involving a measure data, although Liang and Rodrigues \cite{Liang} considered a problem involving measure data both on the domain and on the boundary $\Gamma_2$. Some work has been done by Gallou\"{e}t \cite{Gallouet}, where the non linearity lies on the boundary with measure supported on the domain $\Omega$ and on the boundary $ \Gamma_2$. The MBVP in \cite{Gallouet} posessess a weak solution $u$ in $ W^{1,q}(\Omega),~\forall~ 1<q<\frac{N}{N-1}$ and the trace of $u$  on $\Gamma$ lies in $W^{1-\frac{1}{q},q}(\Gamma)$, $\forall 1\leq q< \frac{N}{N-1}$.\\
In this article	we have considered the following two mixed boundary value problems. The first problem (P1) is
	\begin{eqnarray}\label{q1}
	\begin{split}
	-\Delta u-\lambda^2 u&= h~\text{in}~\Omega,\\
	u&= f_1~\text{on}~\Gamma_1,\\
	\frac{\partial u}{\partial \hat{n}}&= f_2~ \text{on}~\Gamma_2,\\
	\end{split}
	\end{eqnarray}
	and
	\begin{eqnarray}\label{q2}	
\begin{split}
-\Delta u-\lambda^2 u&= 0~\text{in}~\mathbb{R}^N\setminus\bar{\Omega},\\
u&= f_1~\text{on}~\Gamma_1,\\
\frac{\partial u}{\partial \hat{n}}&= f_2~ \text{on}~\Gamma_2,\\
\end{split}
	\end{eqnarray}
where $u$ satisfies the following conditions at infinity, i.e. $|x|\rightarrow\infty.$\\
For $\lambda=0$
\begin{equation}\label{zero}
u(x)=O(|x|^{2-N}).
\end{equation}
For $\lambda\neq0$ (Sommerfeld's radiation condition)
\begin{equation}\label{nonzero}
\begin{split}
 u(x)&=O(|x|^{\frac{1-N}{2}})\\
 \frac{\partial u(x)}{\partial |x|}-i\lambda u(x)&=o(|x|^{\frac{1-N}{2}}).
\end{split}
\end{equation}
The second problem (P2) is
	\begin{eqnarray}\label{lq1}
	\begin{split}
	-\Delta u&= \mu~\text{in}~\Omega,\\
	u&= f_1~\text{on}~\Gamma_1,\\
	\frac{\partial u}{\partial \hat{n}}&= f_2~ \text{on}~\Gamma_2,\\
	\end{split}	
	\end{eqnarray}
	and
	\begin{eqnarray}\label{lq2}
	\begin{split}
	-\Delta u&= 0~\text{in}~\mathbb{R}^N\setminus\bar{\Omega},\\
	u&= f_1~\text{on}~\Gamma_1,\\
	\frac{\partial u}{\partial \hat{n}}&= f_2~ \text{on}~\Gamma_2,
	\end{split}
	\end{eqnarray}
where $u$ satisfies 
\begin{equation}\label{infinity}
|u(x)|+|x||\nabla u(x)|=O(|x|^{2-N}),~\text{as}~|x|\rightarrow\infty.
\end{equation}
Throughout the article $\frac{\partial u}{\partial \hat{n}}$ will represent the normal derivative with respect to the outward unit normal $\hat{n}$  to the boundary, $f_1\in H^{\frac{1}{2}}(\Gamma_1)$, $f_2\in H^{-\frac{1}{2}}(\Gamma_2)$ are boundary data, $h\in \widetilde{H}^{-1}(\Omega)$, $\mu$ will denote a bounded Radon measure, $\lambda\in\mathbb{C}$ with $Im(\lambda)\geq0$ and $\lambda^2$ will be different from the eigenvalues of the Laplacian ($-\Delta$). We will, at some places, refer problem $\eqref{q1},\eqref{lq1}$ as interior problems (IP1), (IP2) respectively and $\eqref{q2},\eqref{lq2}$ as exterior problems (EP1), (EP2) respectively. This work is motivated from the work of Chang \cite{Chang} and Stephan \cite{Stephan} where the authors have used the method of layer potentials to show the uniqueness of solution to the homogeneous mixed bounadry value problem in both interior and exterior domains. Chang \cite{Chang} has shown that for $h=0$ and $\lambda=0$ the solution $u$ belongs to $H^1(\Omega)$ for the interior problem and belongs to $H^1_{loc}(\mathbb{R}^N\setminus\bar{\Omega})$ for the exterior problem. This $u$ also satisfies the following inequality.
 $$\norm{u}_{ H^{1/2}(\Gamma_1)}+\norm{\frac{\partial u}{\partial\hat{n}}}_{ H^{-1/2}(\Gamma_2)}\leq C\{\norm{f_1}_{ H^{1/2}(\Gamma_1)}+\norm{f_2}_{ H^{-1/2}(\Gamma_2)}\}$$
where $C$ is independent of $f_1$, $f_2$ and $h$. The novelty of our work is the consideration of two nonhomogenous mixed boundary value problems and a Radon measure $\mu$ as a nonhomogeneous term in (P2), for which the solution space becomes weaker than the Sobolev space $H^1(\Omega)$.
\section{Preliminary definitions and properties of boundary layer potentials}
We will denote several constants by $C$ which can only depend on $\Omega$, $N$ and independent of the indices of the sequences. The value of $C$ can be different from line to line and sometimes, on the same line.\\  
For $1\leq p\leq\infty$ and $k$ be a nonnegative integer, the Sobolev space $\{u\in L^p(\Omega):D^\gamma u\in L^p(\Omega),~\text{for}~|\gamma|\leq k\}$ will be denoted by $W^{k,p}(\Omega)$ \cite{Evans} and the norm on vectors in $W^{k,p}(\Omega)$ is defined as
\begin{center}
	$\norm{u}_{W^{k,p}(\Omega)}=\sum_{|\gamma|\leq k}\norm{D^\gamma u}_{L^p(\Omega)}$
\end{center}
where $\Omega$ is a domain in $\mathbb{R}^N$. We denote  $W_{loc}^{k,p}(\Omega)$ to be the local Sobolev space such that for any compact $K\subset \Omega$, $u\in W^{k,p}(K)$.  For $0<\alpha<1$, we define the Sobolev space $W^{\alpha,p}(\Omega)$ as
$$W^{\alpha,p}(\Omega)=\{u\in L^p(\Omega):\norm{u}_{W^{\alpha,p}(\Omega)}^p=\norm{u}_{L^p(\Omega)}^p+\int_{\Omega}\int_{\Omega}\frac{|u(x)-u(y)|^p}{|x-y|^{n+p\alpha}}dydx<\infty\}.$$
 Let $\Omega$ be a bounded Lipschitz domain in $\mathbb{R}^N$, $N\geq3$. We now introduce the following Sobolev spaces. For $p=2$, $s\in\mathbb{R}$ and $0\leq\alpha\leq 1$,
\begin{enumerate}
	\item$H^s(\mathbb{R}^N)=\{u:\int_{\mathbb{R}^N}(1+|\xi|^2)^{s/2}\widehat{u}(\xi)e^{i2\pi\xi.x}d\xi \in L^2(\mathbb{R}^N)\},~ \widehat{u}$ is the Fourier transform of $u$. This space is a separable Hilbert space.
	\item $H^s(\Omega)=\{u|_{\Omega}:u\in H^s(\mathbb{R}^N)\}$
	\item  $\widetilde{H}^s(\Omega)=\text{Closure of}~C_0^\infty(\Omega)~\text{in}~H^s(\mathbb{R}^N)$. For further details on these Sobolev spaces one may refer to \cite{Hsiao} Chapter 4.
	\item	$H^\alpha(\Gamma)=\begin{cases}
	\{g|_\Gamma: g\in H^{\frac{1}{2}+\alpha}(\Omega)\},&(0<\alpha\leq1)\\
	L^2(\Gamma)&(\alpha=0),\end{cases}$
	\item $H^{\alpha}(\Gamma_i)= \{g|_{\Gamma_i}:g\in H^{\alpha}(\Gamma)\}$,
	\item $\widetilde{H}^{\alpha}(\Gamma_i)= \{g|_{\Gamma_i}:g\in H^{\alpha}(\Gamma),~\text{supp}(g)\subset\overline{\Gamma}_i\},~i=1,2.$
\end{enumerate}
Let $H^{-\alpha}(\Gamma)$ is the dual space of $H^{\alpha}(\Gamma)$, i.e. $H^{-\alpha}(\Gamma)=(H^{\alpha}(\Gamma))^*$. Equivalently  $H^{-\alpha}(\Gamma_i)=(\widetilde{H}^{\alpha}(\Gamma_i))^*$ and $\widetilde{H}^{-\alpha}(\Gamma_i)=(H^{\alpha}(\Gamma_i))^*$, for $i=1,2$. \\
We denote $\langle.,.\rangle_\Gamma$ as the duality pairing between $H^\alpha(\Gamma)$ and $H^{-\alpha}(\Gamma)$ given by $\langle f,g\rangle_\Gamma=\int_\Gamma f(z)g(z)ds_z$ for any $f\in H^\alpha(\Gamma)$ and $g\in H^{-\alpha}(\Gamma)$. Similarly, $\langle.,.\rangle_{\Gamma_i}$ is the duality pairing between $H^{\alpha}(\Gamma_i)$ and $\widetilde{H}^{-\alpha}(\Gamma_i)$ (or $H^{-\alpha}(\Gamma_i)$ and $\widetilde{H}^{\alpha}(\Gamma_i)$), $i=1,2$.\\
Since $\widetilde H^{\alpha}(\Gamma_i)$, $i=1,2$, is a reflexive space, the operator $$J:\widetilde H^{\alpha}(\Gamma_i)\rightarrow (\widetilde H^{\alpha}(\Gamma_i))^{**}=(H^{-\alpha}(\Gamma_i))^*$$ is a bijection. Hence, for any $f'\in(H^{-\alpha}(\Gamma_i))^*$ there exists a unique $f\in \widetilde H^{\alpha}(\Gamma_i)$ such that $J(f)=f'$. For $g\in H^{-\alpha}(\Gamma_i)$ we define $\langle\langle.,.\rangle\rangle$ by the duality pairing between $H^{-\alpha}(\Gamma_i)$ and $(H^{-\alpha}(\Gamma_i))^*$ such that
$$\langle\langle f',g\rangle\rangle=\langle g,f\rangle_{\Gamma_i}.$$ 
\begin{definition}
	An open set $\Omega\in \mathbb{R}^N$ is said to be a Lipschitz domain if for each $P\in\partial\Omega$ there exist a rectangular coordinate system, $(x,z)$ such that $x\in \mathbb{R}^{n-1}, z\in \mathbb{R}$, a neighborhood $N(P)=N\subset\mathbb{R}^N$ and a function $\varphi:\mathbb{R}^{n-1}\rightarrow\mathbb{R}$ such that
	\begin{enumerate}
		\item 	$|\varphi(x)-\varphi(y)|\leq C|x-y|, ~\forall~x,y\in\mathbb{R}^{n-1}$,
		\item $ N\cap\Omega=\{(x,z):z>\varphi(x)\}\cap N$.
	\end{enumerate}
\end{definition}
\begin{definition}
	The Marcinkiewicz space denoted as ${M}^r(\Omega)$ (or weak $L^r(\Omega)$ space), for every $ 0 < r <\infty$, consists of all measurable functions $g:\Omega\rightarrow\mathbb{R}$ such that 
	$$m(\{x\in \Omega:|g(x)|>b\})\leq \frac{C}{b^r},~b>0,C<\infty,$$
\end{definition}
\noindent where $m$ is the Lebesgue measure. In fact in the case of bounded domain $\Omega$, for any fixed $\bar{r}>0$  we observe ${M}^r(\Omega)\subset {M}^{\bar{r}}(\Omega)$ for  $r\geq \bar{r}$. Furthermore, the embeddings 
\begin{equation}\label{marcinq}
L^r(\Omega)\hookrightarrow {M}^r(\Omega)\hookrightarrow L^{r-\epsilon}(\Omega),
\end{equation}
is continuous for every $1<r<\infty$ and $0<\epsilon<r-1$.
\begin{definition}\label{defn}
	\textbf{(Fredholm operator)} Let $X~ \text{and}~ Y$ are two Banach spaces and $A$ is a bounded linear operator from $X$ to $Y$. Then $A$ is said to be a Fredholm operator if its kernel (ker($A$)) and cokernel (coker($A$)$~=Y/Range(A)$  are finite dimensional. 
\end{definition}
\begin{remark}
	\begin{enumerate}
		\item  The ``Fredholmness" of an operator $A$ ensures that Range$(A)$ is closed.
		\item The index of a Fredholm operator $A$ is given by ind$(A)$=dim(ker$(A)$)-dim(coker$(A)$).
	\end{enumerate}
\end{remark}
\noindent The following two theorems are borrowed from \cite{Driver} which show the relationship between a Fredholm operator and a compact operator. 
\begin{theorem}\label{ft1}
	For a bounded linear operator $A:X\rightarrow Y$, the following two statements are equivalent\\
	1. $A$ is a Fredholm operator.\\
	2. $A$ is an invertible modulo compact operators, i.e. there exist compact operators $C_1,C_2$ and an operator $B$ such that $AB=I+C_1$ and $BA=I+C_2$.  
\end{theorem}
\begin{theorem}\label{ft2}
	If $A$ is a Fredholm opertor then ind$(A)=0$ iff $A=A_1+A_2$, where $A_1$ is an invertible operator and $A_2$ is a compact operator.
\end{theorem}
\begin{definition}
	The space of all finite Radon measures on $\Omega\subset\mathbb{R}^N$, is denoted as $\mathcal{M}(\Omega)$. For $\mu\in \mathcal{M}(\Omega)$ we define $$\norm{\mu}_{\mathcal{M}(\Omega)}=\int_{\Omega}d|\mu|,$$
	which is called the `Total variation' norm.
\end{definition}
\noindent We now define the weak solution of the first problem (P1).
\begin{definition}
	Let $X$ and $Y$ are two test function spaces defined as  $X=\{\varphi\in C^1(\bar{\Omega}):\varphi|_{\Gamma_1}=0\}$ and  $Y=\{\zeta\in C_c^1(\mathbb{R}^N\setminus(\Omega)):\zeta|_{\Gamma_1}=0~ \text{and satisfies}~ \eqref{zero}~\text{and}~\eqref{nonzero}\}$.
	A function $u\in W^{1,1}(\Omega)$ is a weak solution to the problem $\eqref{q1}$ if it satisfies 	
	\begin{equation}
	\int_\Omega \nabla u\cdot\nabla \varphi~ -\int_{\Omega}\lambda^2 u\varphi= \int_\Omega h\varphi +\int_{\Gamma_2}f_2\varphi, ~\forall \varphi\in X.\nonumber
	\end{equation}
	Similarly a function $u\in W^{1,1}_{loc}(\mathbb{R}^N\setminus\bar{\Omega})$ is said to be a weak solution of $\eqref{q2}$ if 
	\begin{equation}
	\int_{\mathbb{R}^N\setminus\bar{\Omega}} \nabla u\cdot\nabla \zeta~ -\int_{\mathbb{R}^N\setminus\bar{\Omega}}\lambda^2 u\zeta=-\int_{\Gamma_2}f_2\zeta, ~\forall \zeta\in Y.\nonumber
	\end{equation}
	\begin{remark}\label{remark 1}
		Hereafter, a subsequence of a sequence will be denoted by the same notation as that of the sequence. Further a solution will always refer to a weak solution. 
	\end{remark}
\end{definition}
\noindent We further we denote $\Phi$ as the fundamental solution of 
Helmholtz equation for $N\geq3$ which satisfies $-\Delta \Phi-\lambda^2 \Phi=\delta$, where $\delta$ is the Dirac distribution and $\Phi$ is
\begin{eqnarray}
\begin{split}
\Phi(x,y)&\begin{cases}
\frac{1}{(N-2)w_N}\frac{1}{|x-y|^{N-2}},&\text{for}~ \lambda=0\\ 
 \frac{e^{i\lambda|x-y|}}{4\pi|x-y|},&\text{for}~ \lambda\neq 0,~ N=3\nonumber\\ \frac{i}{4}\bigg(\frac{{\lambda}}{2\pi(|x-y|)}\bigg)^{\frac{N-2}{2}}H^{(1)}_{\frac{N-2}{2}}(\lambda|x-y|),&\text{for}~\lambda\neq0,~N\geq3
\end{cases}
\end{split}
\end{eqnarray} 
for every $x,y\in\mathbb{R}^N$, $x\neq y$. Here $w_N$ is the measure of the unit sphere in $\mathbb{R}^N$ and $H^{(1)}_m$ denotes the Hankel function of the first kind of order $m$. We next define, boundary layer potentials (single layer and double layer) to solve the homogeneous Helmholtz equation in $\mathbb{R}^N$. Let $g_1\in  H^{\alpha}(\Gamma),g_2\in H^{-\alpha}(\Gamma)$ for some $0\leq\alpha\leq1$, then the single layer potential is given by,
\begin{align}\label{o}
v_1(x)&= \mathbb{S}_{\lambda}g_2(x)\nonumber\\
&=\int_\Gamma g_2(y)\Phi(x-y) dy ,~\forall~ x\in \mathbb{R}^N\setminus\Gamma
\end{align}
and the double layer potential is by
\begin{align}\label{p}
v_2(x)&=\mathbb{K}_{\lambda}g_1(x)\nonumber\\
&= \int_\Gamma g_1(y)\frac{\partial}{\partial \hat{n}_y}\Phi(x-y) dy,~\forall~x\in\mathbb{R}^N\setminus\Gamma.
\end{align}
where $ \hat{n}_y$ denotes the unit outward normal to the boundary $\Gamma$. We can see that for $x\in\mathbb{R}^N\setminus\Gamma$ the above two kernels are $C^\infty$ functions on $\Gamma$.\\
If $P\in\Gamma$, then $X(P)$ denotes a cone with vertex at $P$ such that one component is in $\Omega$ which is denoted by $X_i(P)$ and the other is in $\mathbb{R}^N\setminus\bar{\Omega}$ denoted by $X_e(P)$.
\begin{definition}
	Let $P\in \Gamma$, then we define 
	$$S_{\lambda}g_2(P)=\int_\Gamma g_2(y)\Phi(P-y) dy $$ and $$K_{\lambda}g_1(P)=\int_\Gamma g_1(y)\frac{\partial}{\partial \hat{n}_y}\Phi(P-y) dy.$$
\end{definition}
\noindent According to the Lemma 3.8 of \cite{Costabel} the boundary values of the two potentials in $\eqref{o}$ and $\eqref{p}$ are given by 
\begin{eqnarray}\label{single1}
\begin{split}
v_1^i(P)&=\lim_{X_i(P),x\rightarrow P}\mathbb{S}_{\lambda}g_2(x)\\&=S_{\lambda}g_2(P),
\end{split}\\
\begin{split}\label{single2}
v_1^e(P)&=\lim_{X_e(P),x\rightarrow P}\mathbb{S}_{\lambda}g_2(x)\\&=S_{\lambda}g_2(P)
\end{split}
\end{eqnarray}
and
\begin{equation}\label{double1}
\begin{split}
v_2^i(P)&=\lim_{X_i(P),x\rightarrow P}\mathbb{K}_{\lambda}g_1(x)\\ &=\left(-\frac{1}{2}I+K_{\lambda}\right)g_1(P),
\end{split}
\end{equation}
\begin{equation}\label{double2}
\begin{split}
v_2^e(P)&=\lim_{X_e(P),x\rightarrow P}\mathbb{K}_{\lambda}g_1(x)\\ &=\left(\frac{1}{2}I+K_{\lambda}\right)g_1(P).
\end{split}
\end{equation}
\noindent In case of inhomogeneous Helmholtz equation $-\Delta u-\lambda^2u=h$ in $\Omega$, where $h\in \widetilde{H}^{-1}(\Omega)$. The Newton potential (or Volume potential) appears in the form, \\
$$\mathbb{N}_\lambda h(x)= \int_{\Omega} \Phi(x-y) h(y) dy, ~x\in \mathbb{R}^N.$$
It is well known that the Newton potential $\mathbb{N}_\lambda:\widetilde{H}^{-1}(\Omega)\rightarrow H^{1}(\Omega)$ is a continuous map by \cite{Mclean,Steinbach}. From \cite{Necas} we know the Dirichlet trace operator, $\gamma_D:H^1(\Omega)\rightarrow H^{\frac{1}{2}}(\Gamma)$ and the Neumann trace operator, $\gamma_N:H^1(\Omega)\rightarrow H^{-\frac{1}{2}}(\Gamma)$, are continuous operators.
 The Dirichlet trace operator of $\mathbb{N}_\lambda$ denoted as $\gamma_D\mathbb{N}_\lambda$ is given by
$$\gamma_D\mathbb{N}_\lambda(h(P))=\lim_{x\rightarrow P}\mathbb{N}_\lambda h(x), ~\forall P\in \Gamma.$$ Thus 
\begin{align}
\norm{\gamma_D\mathbb{N}_\lambda(h)}_{H^{\frac{1}{2}}(\Gamma)}&\leq C\norm{\mathbb{N}_\lambda(h)}_{H^1(\Omega)}\nonumber\\
&\leq C \norm{h}_{\widetilde{H}^{-1}(\Omega)}.\nonumber
\end{align}
The Neumann trace of $\mathbb{N}_\lambda$ is denoted as $\gamma_N\mathbb{N}_\lambda$ and hence it satisfies 
$$\norm{\gamma_N\mathbb{N}_\lambda(h)}_{H^{-\frac{1}{2}}(\Gamma)}\leq C \norm{h}_{\widetilde{H}^{-1}(\Omega)}.$$
Let us fix $\alpha=\frac{1}{2}$. Consider the single layer potential $v_1(x)=\mathbb{S}_\lambda{g}_2(x),~ \text{for} ~g_2\in H^{-\frac{1}{2}}(\Gamma) $. Then $v_1$ solves the Helmholtz equation in $\mathbb{R}^N\setminus\Gamma$. Thus $v_1\in H^1(\Omega)$ for (IP1), $v_1\in H_{loc}^1(\mathbb{R}^N\setminus\bar{\Omega})$ for (EP1) and satisfies $\eqref{zero}$-$\eqref{nonzero}$ at infinity. We now define the ouward normal derivative of $v_1$, i.e. $\frac{\partial v_1}{\partial \hat{n}}$ that belongs to $H^{-\frac{1}{2}}(\Gamma)$. Let us choose $h_1,h_2\in H^{\frac{1}{2}}(\Gamma)$. We will denote $h_1^*,~h_2^*$ to be the extensions of $h_1$, $h_2$ respectively such that
\begin{align}\label{cont}
\norm{h_1^*}_{H^1(\Omega)}\leq C 
\norm{h_1}_{H^{\frac{1}{2}}(\Gamma)}, ~
\norm{h_2^*}_{H^1(\mathbb{R}^N\setminus\bar\Omega)}\leq C \norm{h_2}_{H^{\frac{1}{2}}(\Gamma)}
\end{align}
for some constant $C>0$ which does not depend on $h_1$ and $h_2$ by \cite{Jonsson}. Define
\begin{equation}\label{normal}
\begin{split}
\left\langle\frac{\partial v_1}{\partial \hat{n}},h_1\right\rangle_\Gamma&= \int_{\Omega}\nabla v_1\cdot\nabla h_1^*-\int_{\Omega}\lambda^2v_1h_1^*,\\
\left\langle\frac{\partial v_1}{\partial \hat{n}},h_2\right\rangle_\Gamma&= -\int_{\mathbb{R}^N\setminus\bar{\Omega}}\nabla v_1\cdot\nabla h_2^*+\int_{\mathbb{R}^N\setminus\bar{\Omega}}\lambda^2v_1h_2^*,
\end{split}
\end{equation}
We have from Costabel et al. \cite{Costabel} that for every $P\in \Gamma$,
\begin{equation}\label{normal1}
\begin{split}
\frac{\partial v_1^i(P)}{\partial \hat{n}}&=\lim_{X_i(P),x\rightarrow P}\frac{\partial v_1(x)}{\partial \hat{n}}\\&= \bigg(\frac{1}{2}I+K_\lambda^*\bigg)g_2(P)
\end{split}
\end{equation}
and
\begin{equation}\label{normal2}
\begin{split}
 \frac{\partial v_1^e(P)}{\partial \hat{n}}&=\lim_{X_e(P),x\rightarrow P}\frac{\partial v_1(x)}{\partial \hat{n}}\\&=\bigg(-\frac{1}{2}I+K_\lambda^*\bigg)g_2(P)
\end{split} 
\end{equation}
where $K_{\lambda}^*$ is the adjoint operator of $K_{\lambda}$ defined as 
\begin{equation}
K_{\lambda}^*g_1(P)=\int_\Gamma\frac{\partial}{\partial \hat{n}_P}\Phi
(P-y)g_1(y)dy.\nonumber
\end{equation}
Similarly, in case of double layer potential $v_2(x)=\mathbb{K}_\lambda{g_1}(x)$ for ${g_1}\in H^{\frac{1}{2}}(\Gamma)$, we have $\frac{\partial v_2}{\partial\hat{n}}\in H^{-\frac{1}{2}}(\Gamma)$ which satisfies $\eqref{normal}$ and $\eqref{cont}$. Let us define an operator $D_\lambda:H^{\frac{1}{2}}(\Gamma)\rightarrow H^{-\frac{1}{2}}(\Gamma)$ as in \cite{Costabel} such that for every $P\in\Gamma$,
\begin{equation}
D_\lambda{g_1}(P)=\frac{\partial}{\partial \hat{n}_P}\mathbb{K}_\lambda g_1(P)
\end{equation}
and
\begin{align}\label{1double}
 \lim\limits_{X_i(P),x\rightarrow P}\frac{\partial}{\partial \hat{n}_x}\mathbb{K}_\lambda g_1(x)&=\lim\limits_{X_e(P),x\rightarrow P}\frac{\partial}{\partial \hat{n}_x}\mathbb{K}_\lambda g_1(x)\nonumber\\&=D_\lambda{g_1}(P).
\end{align}
\begin{lemma}\label{lemma1} The operators 
\begin{enumerate}
		\item $S_{\lambda}:H^{-\frac{1}{2}}(\Gamma)\rightarrow H^{\frac{1}{2}}(\Gamma)$,
		\item  $\big(\pm\frac{1}{2}I+{K_{\lambda}}\big):H^{\frac{1}{2}}(\Gamma)\rightarrow H^{\frac{1}{2}}(\Gamma)$,
		\item $\big(\pm\frac{1}{2}I+{K^*_{\lambda}}\big):H^{-\frac{1}{2}}(\Gamma)\rightarrow H^{-\frac{1}{2}}(\Gamma)$,
		\item $D_{\lambda}:H^{\frac{1}{2}}(\Gamma)\rightarrow H^{-\frac{1}{2}}(\Gamma)$
	\end{enumerate}
	are continuous by $\cite{Costabel}$.
\end{lemma}
\subsection{Derivation of representation formulae}
Let $\Omega$ be a bounded Lipschitz domain in $\mathbb{R}^N$ and $\overline{\Gamma}_1\cup\Gamma_2=\Gamma_1\cup\overline{\Gamma}_2=\Gamma$, $\Gamma_1\cap\Gamma_2=\emptyset$. For   $g_1\in \widetilde{H}^{\frac{1}{2}}(\Gamma_i),~ i=1,2$, we denote the zero extension function $\widetilde{g}_1$ of $g_1$ by 
\[ \widetilde{g}_1= \begin{cases} 
g_1 & in ~\Gamma_i \\
0 & in ~\Gamma\setminus\Gamma_i,~i=1,2.
\end{cases}
\]
Clearly, $\widetilde{g}_1\in H^{\frac{1}{2}}(\Gamma)$. Similarly, for $g_2\in \widetilde{H}^{-\frac{1}{2}}(\Gamma_i),~ i=1,2$, we extend $g_2$ to a function $\widetilde{g}_2\in {H}^{-\frac{1}{2}}(\Gamma)$. We now introduce the following operators.
$$S_{ij}:\widetilde{H}^{-\frac{1}{2}}(\Gamma_i)\rightarrow H^{\frac{1}{2}}(\Gamma_j), ~ S_{ij}g_2=S_\lambda\widetilde{g}_2|_{\Gamma_j} ~\text{for}~ g_2\in \widetilde{H}^{-\frac{1}{2}}(\Gamma_i),$$
$$K_{ij}:\widetilde{H}^{\frac{1}{2}}(\Gamma_i)\rightarrow H^{\frac{1}{2}}(\Gamma_j), ~ K_{ij}g_1=K_\lambda\widetilde{g}_1|_{\Gamma_j} ~\text{for}~ g_1\in \widetilde{H}^{\frac{1}{2}}(\Gamma_i),$$
$$K^*_{ij}:\widetilde{H}^{-\frac{1}{2}}(\Gamma_i)\rightarrow H^{-\frac{1}{2}}(\Gamma_j), ~ K^*_{ij}g_2=K_\lambda^*\widetilde{g}_2|_{\Gamma_j} ~\text{for}~ g_2\in \widetilde{H}^{-\frac{1}{2}}(\Gamma_i),$$
$$D_{ij}:\widetilde{H}^{\frac{1}{2}}(\Gamma_i)\rightarrow H^{-\frac{1}{2}}(\Gamma_j), ~ D_{ij}g_1=D_\lambda\widetilde{g}_1|_{\Gamma_j} ~\text{for}~ g_1\in \widetilde{H}^{\frac{1}{2}}(\Gamma_i).$$
Let $u\in H^1(\Omega)$ be a solution to the Helmholtz equation   $-\Delta u-\lambda^2u=h$ in $\Omega$ and  $u\in H_{loc}^1(\mathbb{R}^N\setminus\bar{\Omega})$ satisfies  $-\Delta u-\lambda^2u=0$ in $\mathbb{R}^N\setminus\bar{\Omega}$ along with $\eqref{zero}-\eqref{nonzero}$. From the Green's second identity we have
$$\int_{\Omega}u\Delta v-v\Delta u =\int_{\Gamma}u\frac{\partial v}{\partial \hat{n}}-v\frac{\partial u}{\partial \hat{n}}~.$$ 
When we replace $v$ with $\Phi$, the fundamental solution of Helmholtz equation, we obtain the following.
\begin{align}\label{green1}
\int_{\Omega}u(y)\Delta \Phi(x,y)-\Delta u(y) \Phi(x,y)&=\int_{\Gamma}u(y)\frac{\partial \Phi(x,y)}{\partial \hat{n}}-\Phi(x,y)\frac{\partial u(y)}{\partial \hat{n}}\nonumber\\
u(x)&=\int_{\Omega}\Phi(x,y) h(y) -\int_{\Gamma}u(y)\frac{\partial \Phi(x,y)}{\partial \hat{n}}-\Phi
(x,y)\frac{\partial u(y)}{\partial \hat{n}}~.
\end{align}
Let $B_r=\{z\in \mathbb{R}^N: |z|=r\}$ and $D_r=\{x\in\mathbb{R}^N\setminus\bar{\Omega}:|x|<r\}$.
On applying the Green's second identity in the domain $D_r$ we get
\begin{align}\label{green2}
u(x)&=-\int_{D_r}u(y)\Delta \Phi(x,y)-\Delta u(y) \Phi(x,y)\nonumber\\
&=-\int_{B_r}u(y)\frac{\partial \Phi(x,y)}{\partial \hat{n}}-\Phi(x,y)\frac{\partial u(y)}{\partial \hat{n}}+\int_{\Gamma}u(y)\frac{\partial \Phi(x,y)}{\partial \hat{n}}-\Phi(x,y)\frac{\partial u(y)}{\partial\hat{n}}~.
\end{align}
On passing the limit $r\rightarrow \infty$ and by using $\eqref{zero}-\eqref{nonzero}$ we see that $$\int_{B_r}u(y)\frac{\partial \Phi(x,y)}{\partial \hat{n}}-\Phi(x,y)\frac{\partial u(y)}{\partial \hat{n}}\rightarrow 0.$$
Let us denote the Cauchy data as $(\phi,\psi)\in H^{\frac{1}{2}}(\Gamma)\times H^{-\frac{1}{2}}(\Gamma)$, where $u|_\Gamma=\phi$ and ${\frac{\partial u}{\partial \hat{n}}}\big|_\Gamma=\psi$.
On combining $\eqref{green1}$ and $\eqref{green2}$, we can express $u$ as
\begin{equation}\label{represerntation}
\begin{split}
u(x)=\begin{cases}\mathbb{N}_\lambda h(x)-\mathbb{K}_\lambda\phi(x)+\mathbb{S}_\lambda\psi(x),&\text{if}~ x\in\Omega\\
 \mathbb{K}_\lambda\phi(x)-\mathbb{S}_\lambda\psi(x),&\text{if}~ x\in\mathbb{R}^N\setminus\bar\Omega.
\end{cases}
\end{split}
\end{equation}
Consider (P1), with the  boundary data $u|_{\Gamma_1}=f_1$ and $\frac{\partial u}{\partial\hat{n}}|_{\Gamma_2}=f_2$, where $f_1\in H^{\frac{1}{2}}(\Gamma_1), f_2\in H^{-\frac{1}{2}}(\Gamma_2)$. For simplicity, we restrict ourselves to the interior mixed boundary value problem $\eqref{q1}$. Obviously the corresponding results for the exterior problem $\eqref{q2}$ are obtained by only slight modifications. Furthermore, we say $\mathring{f}_1, \mathring{f}_2$ are the extensions of $f_1$ and $f_2$ respectively which satisfies
\begin{equation}\label{e1}
\norm{\mathring{f}_1}_{H^{\frac{1}{2}}(\Gamma)}\leq C \norm{f_1}_{H^{\frac{1}{2}}(\Gamma_1)}
\end{equation}
and
\begin{equation}\label{e2}
\norm{\mathring{f}_2}_{H^{-\frac{1}{2}}(\Gamma)}\leq C\norm{f_2}_{H^{-\frac{1}{2}}(\Gamma_2)}.
\end{equation}
The above extension is possible since we know $\partial\Gamma_1=\partial\Gamma_2$ and $\Gamma$ is Lipschitz [3].
Let us define $\phi=\mathring{f}_1+\widetilde{g}_1$ and $\psi=\mathring{f}_2+\widetilde{g}_2$, where $\widetilde{g}_1$ and $\widetilde{g}_2$ are arbitrary functions in ${H}^{\frac{1}{2}}(\Gamma)$ and  ${H}^{-\frac{1}{2}}(\Gamma)$ respectively. Here  $\widetilde{ g}_1$ is the zero extension
of $g_1\in \widetilde{H}^{\frac{1}{2}}(\Gamma_2)$ and $\widetilde{ g}_2$ is the zero extension
of $g_2\in \widetilde{H}^{-\frac{1}{2}}(\Gamma_1)$. The representation $\eqref{represerntation}$ is used to express the solutions of problem $\eqref{q1}$ as
\begin{equation}\label{interior}
u(x)=\mathbb{N}_\lambda h(x)-\mathbb{K}_\lambda(\mathring{f}_1+\widetilde{g}_1)(x)+\mathbb{S}_\lambda(\mathring{f}_2+\widetilde{g}_2)(x).
\end{equation}
On restricting the equation $\eqref{interior}$ to $\Gamma$  we get, 
$$\mathring{f}_1+\widetilde{g}_1= \gamma_D\mathbb{N}_\lambda h-\left(-\frac{1}{2}I+K_\lambda\right)(\mathring{f}_1+\widetilde{g}_1)-S_\lambda(\mathring{f}_2+\widetilde{g}_2).$$
On $\Gamma_1$ we have the following,
\begin{align}\label{dirichlet}
f_1&= \gamma_D\mathbb{N}_\lambda h|_{\Gamma_1}-\left(-\frac{1}{2}I+K_\lambda\right)(\mathring{f}_1+\widetilde{g}_1)\Big|_{\Gamma_1}+S_\lambda(\mathring{f}_2+\widetilde{g}_2)\big|_{\Gamma_1}\nonumber\\
&= \gamma_D\mathbb{N}_\lambda h|_{\Gamma_1}-K_{21}g_1-\left(-\frac{1}{2}I+K_\lambda\right)\mathring{f}_1\Big|_{\Gamma_1}+S_{11}g_2+S_\lambda\mathring{f}_2\big|_{\Gamma_1}\nonumber\\
K_{21}g_1-S_{11}g_2 &= -f_1+\gamma_D\mathbb{N}_\lambda h|_{\Gamma_1}-\left(-\frac{1}{2}I+K_\lambda\right)\mathring{f}_1\Big|_{\Gamma_1}+S_\lambda\mathring{f}_2\big|_{\Gamma_1}\nonumber\\
&= F^*({f_1},{f_2},h)~(\text{say}).
\end{align}
Taking the Neumann trace of $\eqref{interior}$  we have  $$\mathring{f}_2+\widetilde{g}_2=\gamma_N\mathbb{N}_\lambda h-D_\lambda(\mathring{f}_1+\widetilde{g}_1)+\left(\frac{1}{2}I+K_\lambda^*\right)(\mathring{f}_2+\widetilde{g}_2).$$
Similarly on $\Gamma_2$,
\begin{align}\label{neumann}
f_2&= \gamma_N\mathbb{N}_\lambda h|_{\Gamma_2}-D_\lambda(\mathring{f}_1+\widetilde{g}_1)\big|_{\Gamma_2}+\left(\frac{1}{2}I+K_\lambda^*\right)(\mathring{f}_2+\widetilde{g}_2)\Big|_{\Gamma_2}\nonumber\\
&= \gamma_N\mathbb{N}_\lambda h|_{\Gamma_2}-D_{22}g_1-D\mathring{f}_1\big|_{\Gamma_2}+K^*_{12}g_2+\left(\frac{1}{2}I+K_\lambda^*\right)\mathring{f}_2\Big|_{\Gamma_2}\nonumber\\
D_{22}g_1-K^*_{12}g_2 &=-f_2+\gamma_N\mathbb{N}_\lambda h|_{\Gamma_2}-D_\lambda\mathring{f}_1\big|_{\Gamma_2}+\left(\frac{1}{2}I+K_\lambda^*\right)\mathring{f}_2\Big|_{\Gamma_2}\nonumber \\
&=G^*({f_1},{f_2},h)~(\text{say}).
\end{align}
Clearly $F^*\in H^{\frac{1}{2}}(\Gamma_1)$ and $G^*\in H^{-\frac{1}{2}}(\Gamma_2) $.
Combining equations $\eqref{dirichlet}$ and $\eqref{neumann}$ we get
\[ \left( \begin{array}{cc}
K_{21} & -S_{11} \\
D_{22} & -K^*_{12}
\end{array} \right)
\left( \begin{array}{cc}
g_1 \\
g_2
\end{array} \right)
= \left( \begin{array}{cc}
F^* \\
G^*
\end{array} \right).
\]
We now define a matrix operator $A$ as 
  \[ A=\left( \begin{array}{cc}
K_{21} & -S_{11} \\
D_{22} & -K^*_{12}
\end{array} \right) \]
where, $A:\widetilde{H}^{\frac{1}{2}}(\Gamma_2)\times \widetilde{H}^{-\frac{1}{2}}(\Gamma_1)\rightarrow {H}^{\frac{1}{2}}(\Gamma_1)\times {H}^{-\frac{1}{2}}(\Gamma_2).$
\subsection{Invertibility of layer potentials.}
For the homogeneous Helmholtz equation with $\lambda=0$, the boundary layer operators $S_{0}:H^{-\frac{1}{2}}(\Gamma)\rightarrow H^{\frac{1}{2}}(\Gamma)$ and $\big(-\frac{1}{2}I+{K_{0}}\big):H^{\frac{1}{2}}(\Gamma)\rightarrow H^{\frac{1}{2}}(\Gamma)$ are bijective operators by \cite{Chang}.
\begin{proposition}\label{inver1} This Proposition is from $\cite{Torres}$ which concludes that for $Im(\lambda)>0$
	\begin{enumerate}
		\item	$S_{\lambda}:L^2(\Gamma)\rightarrow H^{1}(\Gamma)$ is invertible.
		\item 	$\big(\pm\frac{1}{2}I+{K_{\lambda}}\big):L^2(\Gamma)\rightarrow L^2(\Gamma)$ is invertible.
		\item $\big(\pm\frac{1}{2}I+{K^*_{\lambda}}\big):L^2(\Gamma)\rightarrow L^2(\Gamma)$ is invertible.
	\end{enumerate}
\end{proposition}
\begin{theorem}\label{inver2}
Let $Im(\lambda)>0$. Then $D_\lambda:H^1(\Gamma)\rightarrow L^2(\Gamma)$ is an invertible operator.
	\begin{proof}
		Let us consider a $g\in L^2(\Gamma)$. From the above Proposition $\ref{inver1}$, $\left(\frac{1}{2}I+K^*_\lambda\right)$ is bijective from  $L^2(\Gamma)$ to itself. Hence, there exists a $g'\in L^2(\Gamma)$ such that $(\frac{1}{2}I+K^*_{\lambda})g'=g$.\\
		Let $v(x)=\mathbb{S}_\lambda\left(-\frac{1}{2}I+K^*_{\lambda}\right)^{-1}g'(x)$. Then $v$ satisfies the homogenous Helmholtz equation in $\mathbb{R}^N\setminus\bar{\Omega}$. Using the properties $\eqref{single2},\eqref{normal2}$ and the decay conditions at infinity $\eqref{zero}-\eqref{nonzero}$  in the exterior domain we have the following representation for $v$. \begin{equation}\label{r}
		v(x)= \mathbb{K}_\lambda f(x)-\mathbb{S}_\lambda g'(x)\nonumber
		\end{equation}
	where, $f={S}_\lambda(-\frac{1}{2}I+K_\lambda^*)^{-1}g'\in H^1(\Gamma)$. Taking the Neumann trace of $v$ we get
		$$ g'=D_\lambda f-\left(-\frac{1}{2}I+K_\lambda^*\right)g'$$
		which implies
		\begin{equation}\label{q9}
	D_\lambda f=\left(\frac{1}{2}I+K_\lambda^*\right)g'. 
		\end{equation}
		Hence, for any $g\in L^2(\Gamma)$, there exists $f\in H^1(\Gamma)$ such that $D_\lambda f=\left(\frac{1}{2}I+K_\lambda^*\right)g'=g.$\\
{\it Claim:} $D_\lambda$ is injective.\\
	Suppose there exists $f\in H^1(\Gamma)$ such that $D_\lambda f=0$ on $\Gamma$. We write $v(x)=\mathbb{K}_\lambda f(x)$, for all  $x$ in $\mathbb{R}^N\setminus\Gamma$. Hence, $v\in H^1(\Omega)$ is a solution of $-\Delta v-\lambda^2v=0$ in $\Omega$ and  $v\in H_{loc}^1(\mathbb{R}^N\setminus\bar{\Omega})$ satisfies  $-\Delta v-\lambda^2v=0$ in $\mathbb{R}^N\setminus\bar{\Omega}$ along with $\eqref{zero}-\eqref{nonzero}$. From the equations $\eqref{double1}$-$\eqref{double2}$ we get
	\begin{align}
	v^i-v^e&=\left(-\frac{1}{2}I+K_\lambda\right)f-\left(\frac{1}{2}I+K_\lambda\right)f\nonumber\\
	&=-f\nonumber
	\end{align}	and from \cite{Torres} we have
	\begin{eqnarray}
	\begin{split}
	D_\lambda f(P)&=\lim_{X_i(P),x\rightarrow P}\frac{\partial}{\partial \hat{n}}v(x)\\
	&=\lim_{X_e(P),x\rightarrow P}\frac{\partial}{\partial \hat{n}}v(x).
\end{split}
	\end{eqnarray}	
	Thus,
	\begin{align}\label{q}
	0&=\langle D_\lambda f,-f\rangle_\Gamma\nonumber\\
	&=\Big\langle D_\lambda f,v^i\Big\rangle_\Gamma-\Big\langle D_\lambda f,v^e\Big\rangle_\Gamma\nonumber\\
	&=\int_{\mathbb{R}^N}|\nabla v|^2-\int_{\mathbb{R}^N}\lambda^2| v|^2
	\end{align}
where the last term in $\eqref{q}$ is due to the fact that $v$ is a solution to the homogeneous Helmholtz equation.	As per our assumption $\lambda^2$ is not an eigen value of ($-\Delta$). Hence, using the conditions $\eqref{zero}-\eqref{nonzero}$ we have $v=0$ a.e. in $\mathbb{R}^N$. Since $v$ is continuous across the boundary, we have $-f=v^i-v^e=0$. This implies
	$f=0~\text{on}~ \Gamma$. So, $D_{\lambda}$ is injective.
	\end{proof}
\end{theorem}
\begin{remark}\label{remark}
The operators $S_\lambda$ and $D_\lambda$ are self-adjoint operators, i.e. $S_\lambda=S_\lambda^*$, $D_\lambda=D_\lambda^*$ (refer Lemma 3.9(a) of $\cite{Costabel})$, where $S_\lambda^*$, $D_\lambda^*$ are the adjoint operators of $S_\lambda$, $D_\lambda$ respectively. Hence, using Proposition $\ref{inver1}$, Theorem $\ref{inver2}$ we obtain $S_\lambda^*:H^{-1}(\Gamma)\rightarrow L^2(\Gamma)$ and $D_\lambda^*:L^2(\Gamma)\rightarrow H^{-1}(\Gamma)$ are invertible operators. Using the properties of real interpolation from Appendix B (Theorem B.2) of $\cite{Mclean}$ on $S_\lambda$, $D_\lambda$ we have
	\begin{enumerate}
		\item $S_\lambda:H^{-\frac{1}{2}}(\Gamma)\rightarrow H^{\frac{1}{2}}(\Gamma)$,
		\item $D_\lambda:H^{\frac{1}{2}}(\Gamma)\rightarrow H^{-\frac{1}{2}}(\Gamma)$
	\end{enumerate}
	are invertible operators.
\end{remark}
\section{Existence and uniqueness results of (P1).}
\begin{theorem}\label{S}
	Let $\Gamma_1\subset\Gamma$, then $S_{11}:\widetilde{H}^{-\frac{1}{2}}(\Gamma_1)\rightarrow {H}^{\frac{1}{2}}(\Gamma_1)$ is a bijective operator.
	\begin{proof}
		We break the proof into three steps.\\
		\textbf{Step 1.} {\it The  operator $S_{11}$ is injective.}\\
		Assume that there exists $g_2\in\widetilde{H}^{-\frac{1}{2}}(\Gamma_1)$ such that $S_{11}g_2=0$ on $\Gamma_1$. We write $v_1(x)=\mathbb{S}_\lambda\widetilde{g}_2(x)$, $\forall x\in \mathbb{R}^N\setminus\Gamma$, where $\widetilde{g}_2\in {H}^{-\frac{1}{2}}(\Gamma)$ is the zero extension of $g_2$. Hence, from the equations $\eqref{single1}$-$\eqref{single2}$ we have $v_1^i=v_1^e\in {H}^{\frac{1}{2}}(\Gamma)$ and from $\eqref{normal1}$, $\frac{\partial v_1^i}{\partial \hat{n}}-\frac{\partial v_1^e}{\partial \hat{n}}=\widetilde{g}_2.$\\
		On replacing $h_1,h_2$ with $v_1^i,v_1^e$ respectively in the equation $\eqref{normal}$ we have
		\begin{align}\label{q5}
			0&= \langle g_2,S_{11}g_2\rangle_{\Gamma_1}\nonumber\\
			&=\langle\widetilde{g}_2,S_\lambda\widetilde{g_2}\rangle_\Gamma\nonumber\\
			&=\Big\langle\frac{\partial v_1^i}{\partial \hat{n}},v_1^i\Big\rangle_\Gamma-\Big\langle\frac{\partial v_1^e}{\partial \hat{n}},v_1^e\Big\rangle_\Gamma\nonumber\\
			&=\int_{\mathbb{R}^N}|\nabla v_1|^2-\int_{\mathbb{R}^N}\lambda^2 |v_1|^2.
		\end{align}
Thus, on using the conditions $\eqref{zero}-\eqref{nonzero}$ we conclude that $v_1=0$ a.e. in $\mathbb{R}^N$. By the continuity of $v_1$ on $\Gamma$, we have  $\widetilde{g}_2=\frac{\partial v_1^i}{\partial \hat{n}}-\frac{\partial v_1^e}{\partial \hat{n}}=0$. This implies $g_2=0$ in $\Gamma_1$ and hence $S_{11}$ is injective.\\
		\textbf{Step 2.} {\it$S_{11}$ is bounded below.}\\
		Suppose there exists a sequence $(g_2^n)\in \widetilde{H}^{-\frac{1}{2}}(\Gamma_1)$ such that $S_{11}g_2^n\rightarrow f$, for some $f\in {H}^{\frac{1}{2}}(\Gamma_1)$. \\
		{\it Case 1.} Assume that $(g_2^n)$ is a bounded sequence in $\widetilde{H}^{-\frac{1}{2}}(\Gamma_1)$. Therefore, there exists a subsequence $(g_2^n)$ and $g_2\in \widetilde{H}^{-\frac{1}{2}}(\Gamma_1)$  such that $(g_2^n)$ converges weakly to $g_2$, i.e.  $g_2^n\overset{w}\rightharpoonup g_2$ in $\widetilde{H}^{-\frac{1}{2}}(\Gamma_1)$. Let $l\in \widetilde{H}^{-\frac{1}{2}}(\Gamma_1)$. Then we have
		\begin{align}
		\langle l,f\rangle_{\Gamma_1}&= \langle l, \lim_{n\rightarrow\infty}S_{11}g_2^n\rangle_{\Gamma_1}\nonumber\\
		&= \lim_{n\rightarrow\infty}\langle l,S_{11}g_2^n\rangle_{\Gamma_1}\nonumber\\ 
		&=\lim_{n\rightarrow\infty} \langle\langle S^*_{11}l,g_2^n\rangle\rangle\nonumber\\
		&=\langle\langle S^*_{11}l, g_2\rangle\rangle\nonumber\\
		&=\langle l, S_{11}g_2\rangle_{\Gamma_1}\nonumber.
		\end{align}
		Since every reflexive space has a unique predual, hence $S_{11}g_2=f$. Therefore, $S_{11}$ has a closed range.\\
		{\it Case 2.} Assume that $(g_2^n)$ is an unbounded sequence in $\widetilde{H}^{-\frac{1}{2}}(\Gamma_1)$. Let us denote
		$$G^n=\frac{g_2^n}{\norm{g_2^n}_{\widetilde{H}^{-\frac{1}{2}}(\Gamma_1)}}.$$ Hence, $\norm{G^n}_{\widetilde{H}^{-\frac{1}{2}}(\Gamma_1)}=1.$
		Therefore, there exists a subsequence $(G^n)$ and $G\in \widetilde{H}^{-\frac{1}{2}}(\Gamma_1)$ such that $G^n\overset{w}\rightharpoonup G$ in $\widetilde{H}^{-\frac{1}{2}}(\Gamma_1)$. Since $S_{11}g_2^n\rightarrow f$ and $\norm{g_2^n}_{\widetilde{H}^{-\frac{1}{2}}(\Gamma_1)}\rightarrow\infty$, we have $S_{11}G^n\rightarrow 0$ in $H^{\frac{1}{2}}(\Gamma_1)$. From  Case 1 it easily follows that $S_{11}G=0$, which further implies $G=0$ by the injectivity of $S_{11}$. Using the invertibility of $S_\lambda$ (refer Remark $\ref{remark}$) we obtain
		\begin{align}\label{contra}
			1&=\norm{G^n}_{\widetilde{H}^{-\frac{1}{2}}(\Gamma_1)}\nonumber\\
			&\leq\norm{\widetilde{G^n}}_{{H}^{-\frac{1}{2}}(\Gamma)}\nonumber\\
			&\leq C\norm{S_{\lambda}(\widetilde{G^n})}_{{H}^{\frac{1}{2}}(\Gamma)}~(\text{for}~ C>0).
		\end{align}
We know that $S_{11}G^n=S_\lambda(\widetilde{G}^n)|_{\Gamma_1}$ and
$S_{12}G^n=S_\lambda(\widetilde{G}^n)|_{\Gamma_2}$.
For $x\neq y$, $\Phi(x-y)$ is a $C^\infty$ function. This implies $S_{12}G^n\rightarrow 0$ in $H^{\frac{1}{2}}(\Gamma_2)$, since $G^n\overset{w}{\rightharpoonup}0$ in $\widetilde{H}^{-\frac{1}{2}}(\Gamma_1)$. Hence, $S_\lambda(\widetilde{G^n})\rightarrow 0$ in $H^{\frac{1}{2}}(\Gamma)$, which is a contradiction to $\eqref{contra}$. Therefore, we conclude that $S_{11}$ has closed range. Thus, $S_{11}$ is bounded below since $S_{11}$ is injective and its range is closed. \\
		\textbf{Step 3.} {\it $S_{11}$ has dense range}.\\
		Assume that $S_{11}^*g_2=0$ for some $g_2\in \widetilde{H}^{-\frac{1}{2}}(\Gamma_1)$. Hence, for any $l\in \widetilde{H}^{-\frac{1}{2}}(\Gamma_1)$,
		\begin{align}
		0&=\langle\langle S_{11}^*g_2,l\rangle\rangle\nonumber\\
		&=\langle g_2, S_{11}l\rangle_{\Gamma_1} .\nonumber
		\end{align}
		Choose $l=g_2$. Then by proceeding on similar lines as in step 1 we get $g_2=0$. Since $\text{Kernel}(S_{11}^*)= {\text{Range}(S_{11})^\perp}=\overline{\text{Range}(S_{11})}^\perp$, the injectivity of $S_{11}^*$ implies $S_{11}$ has dense range.\\
		 Combining the results from the above three steps we conclude that the operator $S_{11}:\widetilde{H}^{-\frac{1}{2}}(\Gamma_1)\rightarrow {H}^{\frac{1}{2}}(\Gamma_1)$ is bijective.
	\end{proof}
\end{theorem}
\begin{theorem}\label{D}
	Let $\Gamma_2\subset\Gamma$, then the operator $D_{22}:\widetilde{H}^{\frac{1}{2}}(\Gamma_2)\rightarrow {H}^{-\frac{1}{2}}(\Gamma_2)$ is invertible.
	\begin{proof}
		Similar to the steps in Theorem $\ref{S}$, we will show that $D_{22}$ is injective and bounded below with a dense range.
		Assume that there exists $g_1\in \widetilde{H}^{\frac{1}{2}}(\Gamma_2)$ such that $D_{22}g_1=0$ on $\Gamma_2$. We now express $v_2(x)=\mathbb{K}_\lambda\widetilde{g}_1(x), ~\forall  x\in\mathbb{R}^N\setminus\Gamma$.  From the equations $\eqref{double1}$ and $\eqref{double2}$ we get
		\begin{align}
			v_2^i-v_2^e&=\left(-\frac{1}{2}I+K_\lambda\right)\widetilde{g}_1-\left(\frac{1}{2}I+K_\lambda\right)\widetilde{g}_1\nonumber\\
			&=-\widetilde{g}_1.\nonumber
		\end{align}		
Thus,
		\begin{align}\label{u}
		0&=\langle D_{22}g_1,-g_1\rangle_{\Gamma_2}\\
		&=\langle D_\lambda\widetilde{g}_1,-\widetilde{g}_1\rangle_\Gamma\nonumber\\
		&=\Big\langle D_\lambda\widetilde{g}_1,v_2^i\Big\rangle_\Gamma-\Big\langle D_\lambda\widetilde{g}_1,v_2^e\Big\rangle_\Gamma~~(\text{from the equation}~\eqref{1double})\nonumber\\
		&=\int_{\mathbb{R}^N}|\nabla v_2|^2-\int_{\mathbb{R}^N}\lambda^2| v_2|^2.\nonumber
		\end{align}
 Hence, using the conditions $\eqref{zero}-\eqref{nonzero}$ we have $v_2=0$ a.e. in $\mathbb{R}^N$, since $\lambda^2$ is not an eigen value of ($-\Delta$). By the continuity of $v_2$ in $x\in\mathbb{R}^N\setminus\Gamma_2$ we have $v_2^i-v_2^e=-\widetilde{g}_1=0$. This implies
		$g_1=0~\text{in}~ \Gamma_2$. So, $D_{22}$ is injective. \\
		On using arguments from Theorem $\ref{S}$, we can show that $D_{22}$ has a closed range and hence it is bounded below. We suppose that $D_{22}^*g_1'=0$ for some $g_1\in\widetilde{H}^{\frac{1}{2}}(\Gamma_2).$ Then for $f\in \widetilde{H}^{\frac{1}{2}}(\Gamma_2)$,
		\begin{align}
		0&=\langle -D_{22}^*g_1',f\rangle_{\Gamma_2}\nonumber\\
		&=\langle\langle- g_1',D_{22}f\rangle\rangle\nonumber\\
		&=\langle -g_1, D_{22}f\rangle_{\Gamma_2}.\nonumber
		\end{align}
	Taking $f=g_1$, then from $\eqref{u}$ we obtain $g_1=0$ in $\Gamma_2$. Hence, $D^*_{22}$ is injective which implies $D_{22}$ has dense range.
		Therefore, $D_{22}$ is an invertible operator.
	\end{proof}
\end{theorem}
\begin{theorem}\label{A}
	The matrix operator $A:\widetilde{H}^{\frac{1}{2}}(\Gamma_2)\times \widetilde{H}^{-\frac{1}{2}}(\Gamma_1)\rightarrow {H}^{\frac{1}{2}}(\Gamma_1)\times{H}^{-\frac{1}{2}}(\Gamma_2)$
	is invertible.
	\begin{proof}
	For any $g_1\in \widetilde{H}^{\frac{1}{2}}(\Gamma_2)$ and $P\in\Gamma_1$, the operator $K_{21}:\widetilde{H}^{\frac{1}{2}}(\Gamma_2)\rightarrow{H}^{\frac{1}{2}}(\Gamma_1)$ is defined as
	\begin{align}\label{s}
	K_{21}g_1(P)&=K_\lambda\widetilde{g}_{1}(P)\nonumber\\
	&=\int_{\Gamma}\frac{\partial}{\partial\hat{n}_y}\Phi(P,y)\widetilde{g}_1(y)dy\nonumber\\
	&=\int_{\Gamma_2}\frac{\partial}{\partial\hat{n}_y}\Phi(P,y){g}_1(y)dy.
	\end{align}
	We can see that the kernel $\frac{\partial\Phi(P,y)}{\partial\hat{n}_y}$ in $\eqref{s}$ is a $C^\infty$ function. Let $(g_1^n)$ be a bounded sequence in $\widetilde{H}^{\frac{1}{2}}(\Gamma_2)$, then there exists a subsequence $(g_1^n)$ and $g_1$ in $\widetilde{H}^{\frac{1}{2}}(\Gamma_2)$ such that $(g_1^n)$ converges weakly to  $g_1$. Hence,
	\begin{align}
	\lim_{n\rightarrow\infty}K_{21}g_1^n(P)&=\lim_{n\rightarrow\infty}\int_{\Gamma_2}\frac{\partial}{\partial\hat{n}_y}\Phi(P,y){g}_1^n(y)dy\nonumber\\
	&=\int_{\Gamma_2}\frac{\partial}{\partial\hat{n}_y}\Phi(P,y) {g}_1(y)dy\nonumber\\
	&=K_{21}g_1(P).\nonumber
	\end{align}
	 Thus, $K_{21}$ is a compact operator. Similarly we can show that the operator $K^*_{12}$ is also compact.\\
		We have
		\begin{align}
		 A&= \left( \begin{array}{cc}
		K_{21} & -S_{11} \\
		D_{22} & -K^*_{12}
		\end{array} \right)\nonumber
	\\	&=\left( \begin{array}{cc}
		K_{21} & 0 \\
		0 & -K_{12}^*
		\end{array} \right)
		+ \left( \begin{array}{cc}
		0 & -S_{11} \\
		D_{22} & 0
		\end{array} \right)\nonumber\\
		&=A_1+A_2\nonumber
			\end{align}
			where $A_1=\left( \begin{array}{cc}
			K_{21} & 0 \\
			0 & -K_{12}^*
			\end{array} \right)$ and $A_2=\left( \begin{array}{cc}
			0 & -S_{11} \\
			D_{22} & 0
			\end{array} \right)$. The matrix $A_2$ is invertible, since $S_{11}$ and $D_{22}$ are invertible operators by Theorem $\ref{S}$ and Theorem $\ref{D}$ respectively. As the operators $D_\lambda$ and $S_\lambda$ are also continuous by \cite{Costabel}, the inverse of $A_2, \text{i.e.}~A_2^{-1}$ is also bounded.  We know the operators $K_{21}~\text{and}~K_{12}^*$ are compact operators and hence $A_1$ is also a compact operator.\\
		Thus, we can write $A_2^{-1}A=A_2^{-1}A_1+I=C_1+I$ and $AA_2^{-1}=A_1A_2^{-1}+I=C_2+I$ where $C_1,C_2$ are compact operators. Using Theorem $\ref{ft1}$, it is equivalent to say that $A$ is a Fredholm operator. This implies ind$(A)=0$ (by Theorem $\ref{ft2}$). Now  to show $A$ is bijective it is sufficient to show $A$ is injective, i.e. dim ker$(A)=0$.   \\
		{\it Claim:} $A$ is injective.\\
		Let us assume that there exist some $g_1\in \widetilde{H}^{\frac{1}{2}}(\Gamma_2)$ and $g_2\in\widetilde{H}^{-\frac{1}{2}}(\Gamma_1)$ such that $A(g_1,g_2)=0.$ Now for $x\in \mathbb{R}^N\setminus\Gamma$, we write
		\begin{eqnarray}
		v(x)=\begin{cases}
		\mathbb{S}_\lambda\widetilde{ g}_2(x)-\mathbb{K}_\lambda\widetilde{ g}_1(x),&\text{if}~x\in \Omega\\
		-\mathbb{S}_\lambda\widetilde{ g}_2(x)+\mathbb{K}_\lambda\widetilde{ g}_1(x),&\text{if}~x\in \mathbb{R}^N\setminus\bar\Omega
		\end{cases}\nonumber
		\end{eqnarray}
		Then $v$ satisfies the following problems
		\begin{eqnarray}\label{q6}
		\begin{split}
		-\Delta v-\lambda^2 v&= 0~\text{in}~\Omega,\\
		v&= 0~\text{on}~\Gamma_1,\\
		\frac{\partial v}{\partial \hat{n}}&= 0~ \text{on}~\Gamma_2,\\
		v &\in H^1(\Omega)
		\end{split}
		\end{eqnarray}
		and
		\begin{eqnarray}\label{q7}
		\begin{split}
		-\Delta v-\lambda^2 v&= 0~\text{in}~\mathbb{R}^N\setminus\bar{\Omega},\\
		v&= 0~\text{on}~\Gamma_1,\\
		\frac{\partial v}{\partial \hat{n}}&= 0~ \text{on}~\Gamma_2,\\
		v &\in H^1_{loc}(\mathbb{R}^N\setminus\bar{\Omega}).
		\end{split}
		\end{eqnarray}
		This implies $v|_\Gamma\in {H}^{\frac{1}{2}}(\Gamma)$, $\frac{\partial v }{\partial \hat{n} }\in {H}^{-\frac{1}{2}}(\Gamma)$ and 
		\begin{align}
		0&=\Big\langle\frac{\partial v}{\partial \hat{n}},v\Big\rangle\nonumber\\
		&=\int_{\Omega}|\nabla v|^2-\int_{\Omega}\lambda^2|v|^2\nonumber
		\end{align}
 and
 \begin{align}
 0&=\Big\langle\frac{\partial v}{\partial \hat{n}},v\Big\rangle\nonumber\\
 &=-\int_{\mathbb{R}^N\setminus\bar{\Omega}}|\nabla v|^2+\int_{\mathbb{R}^N\setminus\bar{\Omega}}\lambda^2|v|^2.\nonumber
 \end{align} 
	Thus, $v=0$ a.e. in $\mathbb{R}^N$,	since $\lambda^2$ is not an eigenvalue of $(-\Delta)$ and $v$ satisfies the radiation conditions at infinity. On $\Gamma_1$, $v=0$ and hence $v^i-v^e=\pm\widetilde{g}_1=0$ and $\frac{\partial v^i}{\partial \hat{n}}-\frac{\partial v^e}{\partial \hat{n}}=\pm\widetilde{g}_2=0$. Thus, $g_1=0$ and $g_2=0$.
	\end{proof}
\end{theorem}
\begin{theorem}\label{main 1}
	The mixed boundary value problem $\eqref{q1}$ with $Im(\lambda)>0$ possesses a unique solution $u$ which is represented as 
	$$u(x)=\mathbb{N}_\lambda h(x)-\mathbb{K}_\lambda(\mathring{f_1}+\widetilde{g_1})(x)+\mathbb{S}_\lambda(\mathring{f_2}+\widetilde{g_2})(x)$$ for unique $g_1\in\widetilde{H}^{\frac{1}{2}}(\Gamma_2)$ and $g_2\in \widetilde{H}^{-\frac{1}{2}}(\Gamma_1)$.
	This solution $u$ also satisfies $\eqref{zero}-\eqref{nonzero}$ and 
	\begin{equation}\label{estimate}
		\norm{u}_{H^{\frac{1}{2}}(\Gamma)}+\norm{\frac{\partial u}{\partial \hat{n}}}_{H^{-\frac{1}{2}}(\Gamma)}\leq C\{\norm{f_1}_{H^{\frac{1}{2}}(\Gamma_1)}+\norm{f_2}_{H^{-\frac{1}{2}}(\Gamma_2)}+\norm{h}_{\widetilde{H}^{-1}(\Omega)}\}.
	\end{equation}
	\begin{proof}
		The solvability  and uniqueness of problem $\eqref{q1}$  depend on the invertibility of the operator $A$. Due to Theorem $\ref{A}$ we know that $A$ is invertible. Hence, there exists a unique pair $(g_1,g_2)\in\widetilde{H}^{\frac{1}{2}}(\Gamma_2)\times \widetilde{H}^{-\frac{1}{2}}(\Gamma_1)$ such that 
		$$A(g_1,g_2)=(F^*,G^*),$$
		where, $F^*\in H^{\frac{1}{2}}(\Gamma_1)$ and $G^*\in H^{-\frac{1}{2}}(\Gamma_2)$ as defined in equations $\eqref{dirichlet}$ and $\eqref{neumann}$. Then we can represent $$u(x)=\mathbb{N}_\lambda h(x)-\mathbb{K}_\lambda\phi(x)+\mathbb{S}_\lambda\psi(x),$$
		where,  $\phi=\mathring{f}_1+\widetilde{g}_1$ and $\psi=\mathring{f}_2+\widetilde{g}_2$ are the Cauchy data for the problem $\eqref{q1}$. Since 
		\begin{align}\label{q8}
			 \left( \begin{array}{cc}
			K_{21} & -S_{11} \\
			D_{22} & -K^*_{12}
			\end{array} \right)
			\left( \begin{array}{cc}
			g_1 \\
			g_2
			\end{array} \right)
			= \left( \begin{array}{cc}
			F^* \\
			G^*
			\end{array} \right),
\end{align}
		so we can write 
		\begin{equation}\label{g1}
		g_1=D_{22}^{-1}(K^*_{12}g_2+G^*).
		\end{equation}
Substituting the value of $g_1$ in the equation $\eqref{q8}$ we get
		\begin{equation}\label{g2}
		(S_{11}-K_{21}D_{22}^{-1}K^*_{12})g_2=K_{21}D_{22}^{-1}G^*-F^*.
		\end{equation}
		We will now show that the operator $H:=S_{11}-K_{21}D_{22}^{-1}K^*_{12}:\widetilde{H}^{-\frac{1}{2}}(\Gamma_1)\rightarrow H^{\frac{1}{2}}(\Gamma_1)$ is invertible. We can then represent $g_1$, $g_2$ as follows. 
	 $$g_1=D_{22}^{-1}G^*+D_{22}^{-1}K_{12}^*H^{-1}(K_{21}D_{22}^{-1}G^*-F^*)$$ and
	 	$$g_2=H^{-1}(K_{21}D_{22}^{-1}G^*-F^*).$$ 
			{\it Claim:} $H=S_{11}-K_{21}D_{22}^{-1}K^*_{12}$ is invertible.\\
			We know $K_{21}$ and $K^*_{12}$ are compact operators. So $K_{21}D_{22}^{-1}K^*_{12}$ is also compact. Using Theorem $\ref{ft2}$ we get ind$(H)=0$, since $S_{11}$ is bijective. Thus we only need to show that $H$ is injective.\\
		Suppose $H(g_2)=0,$ for some $g_2\in \widetilde{H}^{-\frac{1}{2}}(\Gamma_1)$. We express $w(x)=\mathbb{S}_\lambda\widetilde{g_2}(x)-\mathbb{K}_\lambda \widetilde{D_{22}^{-1}K^*_{12}g_2}(x)$. 
			We observe that on $\Gamma_1$,
			\begin{align}
	w& = (S_{11}-K_{21}D_{22}^{-1}K^*_{12})g_2\nonumber\\
				&=H(g_2)\nonumber\\
				&=0\nonumber
			\end{align}
		 and on $\Gamma_2$,
		 \begin{align}
		 \frac{\partial w}{\partial \hat{n}} &= K^*_{12}g_2-D_{22}\{D_{22}^{-1}K^*_{12}\}g_2\nonumber\\
		 &=0.\nonumber
		 \end{align}
 Therefore, $w$ satisfies $\eqref{q6}$ and $\eqref{q7}$. Following the proof of Theorem $\ref{A}$ we conclude that $g_2=0$.  So  $H$ is injective hence invertible.\\
		Furthermore, 
		\begin{align}\label{inq1}
		\norm {u}_{H^{\frac{1}{2}}(\Gamma)}&= \norm{\mathring{f}_1+\widetilde{g}_1}_{H^{\frac{1}{2}}(\Gamma)}\nonumber\\
		& \leq \norm {\mathring{f}_1}_{H^{\frac{1}{2}}(\Gamma)}+\norm{\widetilde{g}_1}_{H^{\frac{1}{2}}(\Gamma)}\nonumber\\
		&\leq C\{\norm{f_1}_{H^{\frac{1}{2}}(\Gamma_1)}+\norm{K^*_{12}g_2+G^*}_{H^{-\frac{1}{2}}(\Gamma_2)}\}~(\text{by}~ \eqref{e1}~ \text{and}~ \eqref{g1})\nonumber\\
		& \leq C\{\norm{f_1}_{H^{\frac{1}{2}}(\Gamma_1)}+\norm{G^*}_{H^{-\frac{1}{2}}(\Gamma_2)}+\norm{g_2}_{H^{-\frac{1}{2}}(\Gamma_1)}\}\nonumber\\
		&\leq C\{\norm{f_1}_{H^{\frac{1}{2}}(\Gamma_1)}+\norm{G^*}_{H^{-\frac{1}{2}}(\Gamma_2)}+\norm{K_{21}D_{22}^{-1}G^*-F^*}_{H^{\frac{1}{2}}(\Gamma_1)}\}~(\text{by}~\eqref{g2})\nonumber\\
		&\leq C\{\norm{f_1}_{H^{\frac{1}{2}}(\Gamma_1)}+\norm{f_2}_{H^{-\frac{1}{2}}(\Gamma_2)}+\norm{h}_{\widetilde{H}^{-1}(\Omega)}\}.~(\text{by}~ \eqref{dirichlet}~\text{and}~ \eqref{neumann})
		\end{align}
	Similarly
	\begin{equation}\label{inq2}
	\norm{\frac{\partial u}{\partial \hat{n}}}_{H^{-\frac{1}{2}}(\Gamma)}\leq C\{\norm{f_1}_{H^{\frac{1}{2}}(\Gamma_1)}+\norm{f_2}_{H^{-\frac{1}{2}}(\Gamma_2)}+\norm{h}_{\widetilde{H}^{-1}(\Omega)}\}.
	\end{equation}
 On combining  inequalities $\eqref{inq1}$ and $\eqref{inq2}$ we get
		$$\norm{u}_{H^{\frac{1}{2}}(\Gamma)}+\norm{\frac{\partial u}{\partial \hat{n}}}_{H^{-\frac{1}{2}}(\Gamma)}\leq C\{\norm{f_1}_{H^{\frac{1}{2}}(\Gamma_1)}+\norm{f_2}_{H^{-\frac{1}{2}}(\Gamma_2)}+\norm{h}_{\widetilde{H}^{-1}(\Omega)}\}.$$	
	\end{proof}
\end{theorem}
\noindent  Proceeding similarly for the exterior problem $\eqref{q2}$ of (P1), the following Theorem can be established.
\begin{theorem}\label{main 2}
For	$Im(\lambda)>0$, problem (P1) with given boundary data $f_1\in H^{\frac{1}{2}}(\Gamma_1)$ and $f_2\in H^{-\frac{1}{2}}(\Gamma_2)$ has a unique solution $u$ which is represented as
		\begin{eqnarray}
		u(x)=\begin{cases}
		\mathbb{N}_\lambda h(x)-\mathbb{K}_\lambda(\mathring{f}_1+\widetilde{g}_1)(x)+\mathbb{S}_\lambda(\mathring{f}_2+\widetilde{g}_2)(x),&\text{if}~x\in \Omega\\
		\mathbb{K}_\lambda(\mathring{f}_1+\widetilde{g}_1)(x)-\mathbb{S}_\lambda(\mathring{f}_2+\widetilde{g}_2)(x),&\text{if}~x\in \mathbb{R}^N\setminus\bar\Omega
		\end{cases}\nonumber
		\end{eqnarray}
 for unique $g_1\in\widetilde{H}^{\frac{1}{2}}(\Gamma_2)$ and $g_2\in \widetilde{H}^{-\frac{1}{2}}(\Gamma_1)$. The solution $u$ belongs to $H^1(\Omega)$ for (IP1) and belongs to $H^1_{loc}(\mathbb{R}^N\setminus\bar{\Omega})$ for (EP1) satisfying the conditions $\eqref{zero}-\eqref{nonzero}$ at infinity. Furthermore, $u$ satisfies $\eqref{estimate}$.
\end{theorem}
 \section{Existence results of (P2).}
 Problem (P2) is a mixed boundary value problem of Poisson equation where $\mu\in\mathcal{M}(\bar\Omega)$ is a bounded Radon measure supported on $\Omega$ and the boundary data are $f_1\in H^{\frac{1}{2}}(\Gamma_1)$ and $f_2\in H^{-\frac{1}{2}}(\Gamma_2)$. 
\begin{definition}\label{weak solution}
	We say a function $u\in W^{1,1}(\Omega)$ is a weak solution to the problem $\eqref{lq1}$ if	
		\begin{equation}
		\int_\Omega \nabla u\cdot\nabla \varphi~ = \int_\Omega \varphi d\mu +\int_{\Gamma_2}f_2\varphi, ~\forall \varphi\in X\nonumber
		\end{equation}
		where, $X=\{\varphi\in C^1(\bar{\Omega}):\varphi|_{\Gamma_1}=0\}$ is the test function space. Similarly a function $u\in W^{1,1}_{loc}(\mathbb{R}^N\setminus\bar{\Omega})$ is said to be a weak solution of $\eqref{lq2}$ if 
	\begin{equation}
	\int_{\mathbb{R}^N\setminus\bar{\Omega}} \nabla u\cdot\nabla \varphi=-\int_{\Gamma_2}f_2\varphi, ~\forall \varphi\in Z\nonumber
	\end{equation}
	where,  $Z=\{\zeta\in C_c^1(\mathbb{R}^N\setminus\Omega):\zeta|_{\Gamma_1}=0~ \text{and satisfies}~ \eqref{infinity}\}$.
\end{definition}
\noindent We will now approximate $\mu\in\mathcal{M}(\bar\Omega)$ by a smooth sequence $(\mu_n)\subset L^\infty(\Omega)$, in the weak* topology, i.e. 	
\begin{equation}\label{convergence}
\int_{\Omega} g~ d\mu_n\rightarrow \int_{\Omega} g ~d\mu, ~\forall g\in C(\bar\Omega).
\end{equation} 
 In order to show the existence of solutions to (P2), we consider the \textquoteleft {\it approximating}\textquoteright ~problems to $\eqref{lq1}$-$\eqref{lq2}$ which are as follows.
\begin{eqnarray}\label{q3}
\begin{split}
-\Delta u_n&= \mu_n~\text{in}~\Omega,\\
u_n&= f_1~\text{on}~\Gamma_1,\\
\frac{\partial u_n}{\partial \hat{n}}&= f_2~ \text{on}~\Gamma_2,\\
\end{split}
\end{eqnarray}
and
\begin{eqnarray}\label{q4}
\begin{split}
-\Delta u_n&= 0~\text{in}~\mathbb{R}^N\setminus\bar{\Omega},\\
u_n&= f_1~\text{on}~\Gamma_1,\\
\frac{\partial u_n}{\partial \hat{n}}&= f_2~ \text{on}~\Gamma_2,\\
\end{split}
\end{eqnarray}
These \textquoteleft {\it approximating}\textquoteright~ problems are special cases of (P1) with $\lambda=0$. The weak formulation to \eqref{q3} is 
\begin{equation}\label{formulation}
\int_\Omega \nabla u_n\cdot\nabla \varphi~ = \int_\Omega \varphi \mu_n +\int_{\Gamma_2}f_2\varphi, ~\forall \varphi\in X.
\end{equation}
\begin{theorem}\label{t1}
	The problems $\eqref{q3}$ and $\eqref{q4}$ admit a unique  solution $u_n$ which is represented as
		\begin{eqnarray}\label{t}
	u_n(x)=\begin{cases}
	\mathbb{N}_0 \mu_n(x)-\mathbb{K}_0(\mathring{f}_1+\widetilde{g}_1)(x)+\mathbb{S}_0(\mathring{f}_2+\widetilde{g}_2)(x),&\text{if}~x\in \Omega\\
	\mathbb{K}_0(\mathring{f}_1+\widetilde{g}_1)(x)-\mathbb{S}_0(\mathring{f}_2+\widetilde{g}_2)(x),&\text{if}~x\in \mathbb{R}^N\setminus\bar\Omega
	\end{cases}
	\end{eqnarray}
for a unique pair $(g_1,g_2)\in\widetilde{H}^{\frac{1}{2}}(\Gamma_2)\times \widetilde{H}^{-\frac{1}{2}}(\Gamma_1)$. The solution $u_n$ belongs to $H^1(\Omega)$ for the problem $\eqref{q3}$ and belongs to $H^1_{loc}(\mathbb{R}^N\setminus\bar{\Omega})$ for the problem $\eqref{q4}$ satisfying the radiation condition $\eqref{infinity}$ at infinity. 
\begin{proof}
The invertibility of the operators $S_{11}$ and $D_{22}$ follows from Theorem 3.1 and Theorem 3.2 of \cite{Chang} respectively. Further, following the proof of Theorem $\ref{A}$, Theorem $\ref{main 1}$ one can see that the matrix operator $A$ is invertible and the problems $\eqref{q3}$-$\eqref{q4}$ have a unique solution denoted as $u_n$. The solution $u_n$ can be represented as in $\eqref{t}$ by Theorem $\ref{main 2}$ and satisfies the condition $\eqref{infinity}$ at infinity.
 	\end{proof}
\end{theorem}
\noindent Now to show that problem in $\eqref{lq1}$, involving measure, possesses a solution $u$ we need to pass the limit $n\rightarrow\infty$ in the weak formulation $\eqref{formulation}$.
 \begin{lemma}\label{bdd}
 	Let us suppose that  $u_n$ is a solution of problem $\eqref{q3}$ with $f_1\in H^{\frac{1}{2}}(\Gamma_1)~\text{and}~ f_2\in H^{-\frac{1}{2}}(\Gamma_2)$. Then the sequence $(u_n)$ is bounded in $W^{1,q}(\Omega)~\forall q<\frac{N}{N-1}$.
 	\begin{proof}
 		From the continuous embedding $\eqref{marcinq}$ we have
 		\begin{equation}\label{marcink}
 		L^{\frac{N}{N-1}}(\Omega)\hookrightarrow {M}^{\frac{N}{N-1}}(\Omega)\hookrightarrow L^{\frac{N}{N-1}-\epsilon}(\Omega).
 		\end{equation}
 		If we can show that $(u_n)$ is bounded in ${M}^{\frac{N}{N-1}}(\Omega)$, then this will also imply $(u_n)$ to be bounded in $L^{\frac{N}{N-1}-\epsilon}(\Omega)$. More precisely, we can say $(u_n)$ to be bounded in $L^q(\Omega)$, for every $q<\frac{N}{N-1}$.\\
 		{\it Claim:} The sequences $(u_n)$ and $(\nabla u_n)$ are bounded in ${M}^{\frac{N}{N-1}}(\Omega)$.\\
 	Let us fix a constant  $a>0$. The truncation function of $u_n$ is defined as
 	$$T_a(u_n)=\max\{-a,\min\{a,u_n\}\}.$$
 		Choose $\varphi\in H^1(\Omega)$. Then from Theorem $\ref{t1}$ we have $\frac{\partial u_n}{\partial\hat{n}}\Big|_{\Gamma}=\mathring{f_2}+\widetilde{g_2}$, for a unique ${g_2}\in \widetilde{H}^{-\frac{1}{2}}(\Gamma_1)$. The weak formulation of $\eqref{q3}$  becomes
 		$$\int_{\Omega}\nabla u_n.\nabla\varphi=\int_{\Omega}\varphi\mu_n+\int_{\Gamma}(\mathring{f_2}+\widetilde{g_2})\varphi$$
 		Consider $\varphi =T_a(u_n)$, then we get
 		\begin{align}\label{truncation}
 		\int_{\Omega}|\nabla T_a(u_n)|^2&\leq
 		 \int_{\Omega}\nabla u_n.\nabla (T_a(u_n))\nonumber\\ &=\int_{\Omega}T_a(u_n)\mu_n+\int_{\Gamma}(\mathring{f_2}+\widetilde{g_2})T_a(u_n)\nonumber\\
 		&\leq a\int_{\Omega}\mu_n+a\int_{\Gamma}(\mathring{f_2}+\widetilde{g_2}) \nonumber\\
 		& \leq Ca.
 		\end{align}
 		Consider 
 		\begin{align*}
 		\{|\nabla u_n|\geq  k\} & = \{|\nabla u_n|\geq  k,u_n< a\} \cup \{|\nabla u_n| \geq k,u_n \geq  a\}
 		\\& \subset \{|\nabla u_n|\geq  k,u_n <a\} \cup \{u_n \geq a\} \\&= B_1\cup B_2\subset \Omega
 		\end{align*}
 	where $B_1= \{|\nabla u_n|\geq  k,u_n <a\}$ and $B_2=\{u_n \geq a\}$. Hence, due to the subadditivity property of the Lebesgue measure \textquoteleft $m$\textquoteright ~we have
 		\begin{equation}\label{sub}
 		m(\{|\nabla u_n|\geq k\}) \leq m(B_1) + m(B_2). 
 		\end{equation}
 		Using the Sobolev inequality, we have
 		\begin{align}\label{sobolev}
 		\left(\int_\Omega |T_a(u_n)|^{2^*}\right)^{\frac{2}{2^*}}&\leq \frac{1}{\lambda_1}\int_{\Omega}|\nabla T_a(u_n)|^2 \nonumber\\
 		&\leq Ca
 		\end{align}
 		where $\lambda_1$ is the first eigenvalue of ($-\Delta$). Now we restrict the above inequality $\eqref{sobolev}$ on $B_2$ to get 
 		\begin{align}
 		& a^2m(\{u_n\geq a\})^{\frac{2}{2^*}}\leq Ca,~(\text{Since}~ T_a(u_n)=a ~\text{in}~ B_2).\nonumber
 		\end{align}
 		Thus, 
 		\begin{equation}
 		m(\{u_n\geq a\})\leq \frac{C}{a^\frac{N}{N-2}},~ \forall a\geq1.\nonumber
 		\end{equation}
 		Hence, $(u_n)$ is bounded in ${M}^{\frac{N}{N-2}}(\Omega)$ and also bounded in ${M}^{\frac{N}{N-1}}(\Omega)$. Similarly on restricting $\eqref{truncation}$ on $B_1$, we have
 		\begin{align}
 			m(\{|\nabla u_n|\geq  k,u_n< a\})&\leq \frac{1}{k^2}\int_\Omega |\nabla T_a(u_n)|^2\nonumber\\
 			&\leq \frac{Ca}{k^2}, ~\forall a>1.\nonumber
 		\end{align}
 		Now the inequality $\eqref{sub}$ becomes
 		\begin{align}
 			m(	\{|\nabla u_n|\geq k\}) &\leq m(\{B_1\}) + m(B_2)\nonumber\\&\leq \frac{Ca}{k^2} + \frac{C}{a^\frac{N}{N-2}},~ \forall a>1.\nonumber
 		\end{align}
 		On choosing $a=k^{\frac{N-2}{N-1}}$ we get
 		$$ m(\{|\nabla u_n|\geq k\})\leq \frac{C}{k^\frac{N}{N-1}},~ \forall k\geq 1.$$
 		So $(\nabla u_n)$ is bounded in ${M}^{\frac{N}{N-1}}(\Omega)$. Hence, we conclude that $(u_n)$ is bounded in $W^{1,q}(\Omega)$ for every $q<\frac{N}{N-1}$.
 	\end{proof}
 \end{lemma}
 \begin{theorem}\label{main3}
 	There exists a weak solution $u$ of $\eqref{lq1}$ in $W^{1,q}(\Omega), ~\forall q<\frac{N}{N-1}$.
 \begin{proof}
 	According to Lemma $\ref{bdd}$, $(u_n)$ is bounded in $W^{1,q}(\Omega)$ which is a reflexive space. This implies that there exists a function $u\in W^{1,q}(\Omega)$ such that $u_n$ converges weakly to $u$, i.e. $u_n\overset{w}{\rightharpoonup} u$ in $W^{1,q}(\Omega), ~\forall q<\frac{N}{N-1}$.\\
 	Thus for $\varphi\in X$,
 	$$\lim_{n\rightarrow +\infty} \int_{\Omega} \nabla u_n . \nabla\varphi = \int_{\Omega}\nabla u .\nabla\varphi.$$
 	The sequence $(\mu_n)$ converges to $\mu$ in the weak* topology in the sense given in $\eqref{convergence}$.
 	On passing the limit $n\rightarrow\infty$ in the weak formulation  $\eqref{formulation}$ involving $\mu_n$ we obtain 
 	\begin{equation}
 	\int_{\Omega}\nabla u.\nabla\varphi=\int_{\Omega}\varphi d\mu+\int_{\Gamma_2}f_2\varphi, ~\forall\varphi\in X.\nonumber
 	\end{equation}
 	Hence, a weak solution of $\eqref{lq1}$ in $W^{1,q}(\Omega)$ for every $q<\frac{N}{N-1}$ is guaranteed. 
 \end{proof}
\end{theorem}
\noindent Now with the consideration of Theorem $\ref{t1}$ and Theorem $\ref{main3}$ we state our main result which is as follows. 
 \begin{theorem}
 	There exists a weak solution $u$ of (P2) with $\mu\in\mathcal{M}(\bar\Omega)$, with support is in $\Omega$, as a nonhomogeneous term,  $f_1\in H^{\frac{1}{2}}(\Gamma_1)$ and  $f_2\in H^{-\frac{1}{2}}(\Gamma_2)$. The solution $u$ belongs to $W^{1,q}(\Omega),~\forall q<\frac{N}{N-1}$ for (IP2)  and belongs to $H^1_{loc}(\mathbb{R}^N\setminus\bar{\Omega})$ for (EP2) satisfying equation $\eqref{infinity}$ at infinity. 
 \end{theorem}
 \section*{Acknowledgement}
 The author Akasmika Panda thanks the financial assistantship received from the Ministry of Human
 Resource Development (M.H.R.D.), Govt. of India. Both the authors also acknowledge the facilities received from the Department of mathematics, National Institute of Technology Rourkela.
 
\end{document}